# GENERALIZATIONS OF THE FAMILYWISE ERROR RATE

By E. L. Lehmann and Joseph P. Romano

*University of California, Berkeley and Stanford University*

Consider the problem of simultaneously testing null hypotheses $H_1, \ldots, H_s$. The usual approach to dealing with the multiplicity problem is to restrict attention to procedures that control the familywise error rate (FWER), the probability of even one false rejection. In many applications, particularly if $s$ is large, one might be willing to tolerate more than one false rejection provided the number of such cases is controlled, thereby increasing the ability of the procedure to detect false null hypotheses. This suggests replacing control of the FWER by controlling the probability of $k$ or more false rejections, which we call the $k$-FWER. We derive both single-step and stepdown procedures that control the $k$-FWER, without making *any* assumptions concerning the dependence structure of the $p$-values of the individual tests. In particular, we derive a stepdown procedure that is quite simple to apply, and prove that it cannot be improved without violation of control of the $k$-FWER. We also consider the false discovery proportion (FDP) defined by the number of false rejections divided by the total number of rejections (defined to be 0 if there are no rejections). The false discovery rate proposed by Benjamini and Hochberg [*J. Roy. Statist. Soc. Ser. B* **57** (1995) 289–300] controls $E(FDP)$. Here, we construct methods such that, for any $\gamma$ and $\alpha$, $P\{FDP > \gamma\} \leq \alpha$. Two stepdown methods are proposed. The first holds under mild conditions on the dependence structure of $p$-values, while the second is more conservative but holds without any dependence assumptions.

**1. Introduction.** In this paper, we will consider the general problem of simultaneously testing a finite number of null hypotheses $H_i$, $i = 1, \ldots, s$. We shall assume that tests for the individual hypotheses are available and the problem is how to combine them into a simultaneous test procedure. The easiest approach is to disregard the multiplicity and simply test each hypothesis at level $\alpha$. However, with such a procedure the probability of one or more false rejections increases with $s$. When the number of true hypotheses is large, we shall be nearly certain to reject some of them.

A classical approach to dealing with this problem is to restrict attention to procedures that control the probability of one or more false rejections. This probability is called the familywise error rate (FWER). Here the term "family" refers to the collection of hypotheses $H_1, \ldots, H_s$ that is being considered for joint testing. Which tests are to be treated jointly as a family depends on the situation.

Once the family has been defined, control of the FWER (at joint level $\alpha$) requires that

$$\text{FWER} \leq \alpha \tag{1}$$

for all possible constellations of true and false hypotheses. A quite broad treatment of methods that control the FWER is presented in [4].

Safeguards against false rejections are of course not the only concern of multiple testing procedures. Corresponding to the power of a single test, one must also consider the ability of a procedure to detect departures from the hypothesis when they do occur. When the number of tests is in the tens or hundreds of thousands, control of the FWER at conventional levels becomes so stringent that individual departures from the hypothesis have little chance of being detected. For this reason, we shall consider an alternative to the FWER that controls false rejections less severely and consequently provides better power.

Specifically, we shall consider the $k$-FWER, the probability of rejecting at least $k$ true null hypotheses. Such an error rate with $k > 1$ is appropriate when one is willing to tolerate one or more false rejections, provided the number of false rejections is controlled.

More formally, suppose data $X$ is available from some model $P \in \Omega$. A general hypothesis $H$ can be viewed as a subset $\omega$ of $\Omega$. For testing $H_i : P \in \omega_i$, $i = 1, \ldots, s$, let $I(P)$ denote the set of true null hypotheses when $P$ is the true probability distribution; that is, $i \in I(P)$ if and only if $P \in \omega_i$. Then, the $k$-FWER, which depends on $P$, is defined to be

$$k\text{-FWER} = P\{\text{reject at least } k \text{ hypotheses } H_i \text{ with } i \in I(P)\}. \tag{2}$$

Control of the $k$-FWER requires that $k$-FWER $\leq \alpha$ for all $P$, that is,

(3) $P\{\text{reject at least } k \text{ hypotheses } H_i \text{ with } i \in I(P)\} \leq \alpha \quad$ for all $P$.

Evidently, the case $k = 1$ reduces to control of the usual FWER.

We will also consider control of the *false discovery proportion* (FDP), defined as the total number of false rejections divided by the total number of rejections (and equal to 0 if there are no rejections). Given a user specified value $\gamma \in (0,1)$, the measure of error control we wish to control is $P\{FDP > \gamma\}$ and we derive methods where this is bounded by $\alpha$.

Recently, there has been a flurry of activity in finding methods that control error rates that are less stringent than the FWER, which is no doubt inspired by the FDR controlling method of Benjamini and Hochberg [1] and applications such as genomic studies where $s$ is so large that control of the FWER is too stringent. For example, Genovese and Wasserman [3] study asymptotic procedures that control the FDP (and the FDR) in the framework of a random effects mixture model. These ideas are extended in [9], where in the context of random fields the number of null hypotheses is uncountable. Korn, Troendle, McShane and Simon [8] provide methods that control both the $k$-FWER and FDP; they provide some justification for their methods, but they are limited to a multivariate permutation model. Alternative methods of control of the $k$-FWER and FDP are given in van der Laan, Dudoit and Pollard [13]; they include both finite sample and asymptotic results. Surprisingly, the methods presented here are distinct from the above techniques. Our methods are not asymptotic and hold under either mild or no assumptions, as long as $p$-values are available for testing each individual hypothesis.

Before describing methods that provide control of the $k$-FWER and FDP, we first recall the notion of a $p$-value, since multiple testing methods are often described by the $p$-values of the individual tests. Consider a single null hypothesis $H : P \in \omega$. Assume a family of tests of $H$, indexed by $\alpha$, with level $\alpha$ rejection regions $S_\alpha$ satisfying

$$(4) \qquad P\{X \in S_\alpha\} \leq \alpha \qquad \text{for all } 0 < \alpha < 1, P \in \omega,$$

and

$$(5) \qquad S_\alpha \subset S_{\alpha'} \qquad \text{whenever } \alpha < \alpha'.$$

Then the $p$-value is defined by

$$(6) \qquad \hat{p} = \hat{p}(X) = \inf\{\alpha : X \in S_\alpha\}.$$

The important property of a $p$-value that will be used later is the following.

LEMMA 1.1. *Assume $\hat{p}$ is defined as above.*

(i) *If $P \in \omega$, then*

$$(7) \qquad P\{\hat{p} \leq u\} \leq u.$$

(ii) *Furthermore,*

$$(8) \qquad P\{\hat{p} \leq u\} \geq P\{X \in S_u\}.$$

*Therefore, if the $S_\alpha$ are such that equality holds in* (4), *then $\hat{p}$ is uniformly distributed on $(0,1)$ when $P \in \omega$.*

PROOF. Assume $P \in \omega$. To prove (i), note that the event $\{\hat{p} \leq u\}$ implies $\{X \in S_{u+\varepsilon}\}$ for any small $\varepsilon > 0$. Therefore,

$$P\{\hat{p} \leq u\} \leq P\{X \in S_{u+\varepsilon}\} \leq u + \varepsilon$$

by assumption (4). Now let $\varepsilon \to 0$. To prove (ii), the event $\{X \in S_u\}$ implies $\{\hat{p} \leq u\}$, and so (8) follows. □

Two classic procedures that control the FWER are the Bonferroni procedure and the Holm procedure. The Bonferroni procedure rejects $H_i$ if its corresponding $p$-value satisfies $\hat{p}_i \leq \alpha/s$. Assuming $\hat{p}_i$ satisfies

(9) $\quad P\{\hat{p}_i \leq u\} \leq u \quad$ for any $u \in (0,1)$ and any $P \in \omega_i$,

the Bonferroni procedure provides strong control of the FWER. Unfortunately, the ability of the Bonferroni procedure to detect cases in which $H_i$ is false will typically be very low since $H_i$ is tested at level $\alpha/s$ which—particularly if $s$ is large—is orders smaller than the conventional $\alpha$ levels.

For this reason procedures are prized for which the levels of the individual tests are increased over $\alpha/s$ without an increase in the FWER. It turns out that such a procedure due to Holm [5] is available under the present minimal assumptions.

The Holm procedure can conveniently be stated in terms of the $p$-values $\hat{p}_1, \ldots, \hat{p}_s$ of the $s$ individual tests. Let the ordered $p$-values be denoted by $\hat{p}_{(1)} \leq \cdots \leq \hat{p}_{(s)}$, and the associated hypotheses by $H_{(1)}, \ldots, H_{(s)}$. Then the Holm procedure is defined stepwise as follows:

*Step* 0. Let $k = 0$.

*Step* 1. If $\hat{p}_{(k+1)} > \alpha/(s-k)$, go to step 2. Otherwise set $k = k+1$ and repeat step 1.

*Step* 2. Reject $H_{(j)}$ for $j \leq k$ and accept $H_{(j)}$ for $j > k$.

The Bonferroni method is an example of a *single-step* procedure, meaning any null hypothesis is rejected if its corresponding $p$-value is less than or equal to a common cutoff value (which in the Bonferroni case is $\alpha/s$). The Holm procedure is a special case of a class of *stepdown* procedures, which we now briefly describe. Let

(10) $\quad\quad\quad\quad\quad\quad \alpha_1 \leq \alpha_2 \leq \cdots \leq \alpha_s$

be constants. If $\hat{p}_{(1)} > \alpha_1$, reject no null hypotheses. Otherwise, if

(11) $\quad\quad\quad\quad\quad\quad \hat{p}_{(1)} \leq \alpha_1, \ldots, \hat{p}_{(r)} \leq \alpha_r,$

reject hypotheses $H_{(1)}, \ldots, H_{(r)}$ where the largest $r$ satisfying (11) is used. That is, a stepdown procedure starts with the most significant $p$-value and continues rejecting hypotheses as long as their corresponding $p$-values are small. The Holm procedure uses $\alpha_i = \alpha/(s - i + 1)$.

**2. Control of the $k$-FWER.** The usual Bonferroni procedure compares each $p$-value $\hat{p}_i$ with $\alpha/s$. Control of the $k$-FWER allows one to increase $\alpha/s$ to $k\alpha/s$, and thereby greatly increase the ability to detect false hypotheses. That such a simple modification results in control of the $k$-FWER is seen in the following result.

THEOREM 2.1. *For testing $H_i: P \in \omega_i$, $i = 1, \ldots, s$, suppose $\hat{p}_i$ satisfies (9). Consider the procedure that rejects any $H_i$ for which $\hat{p}_i \leq k\alpha/s$.*

(i) *This procedure controls the $k$-FWER, so that (3) holds. Equivalently, if each of the hypotheses is tested at level $k\alpha/s$, then the $k$-FWER is controlled.*

(ii) *For this procedure, the inequality (3) is sharp in the sense that there exists a joint distribution for $(\hat{p}_1, \ldots, \hat{p}_s)$ for which equality is attained in (3).*

PROOF. (i) Fix any $P$ and suppose $H_i$ with $i \in I = I(P)$ are true and the remainder false, with $|I|$ denoting the cardinality of $I$. Let $N$ be the number of false rejections. Then, by Markov's inequality,

$$P\{N \geq k\} \leq \frac{E(N)}{k} = \frac{E[\sum_{i \in I(P)} I\{\hat{p}_i \leq k\alpha/s\}]}{k} = \sum_{i \in I(P)} \frac{P\{\hat{p}_i \leq k\alpha/s\}}{k}$$

$$\leq \sum_{i \in I(P)} \frac{k\alpha/s}{k} = |I(P)|\frac{\alpha}{s} \leq \alpha.$$

To prove (ii), consider the following construction. Pick $k$ indices at random without replacement from $\{1, \ldots, s\}$. Call them $J$. Given $i \in J$, let $\hat{p}_i = U_1$, where $U_1$ is uniform on $(0, k/s)$, that is, $U_1 \sim U(0, k/s)$. Given $i \notin J$, let $\hat{p}_i = U_2$, where $U_2$ is independent of $U_1$ and $U_2 \sim U(k/s, 1)$. Then, unconditionally,

$$\hat{p}_i \sim \frac{k}{s}U\left(0, \frac{k}{s}\right) + \left(1 - \frac{k}{s}\right)U\left(\frac{k}{s}, 1\right) \sim U(0, 1).$$

Indeed, if $u \leq k/s$,

$$P\{\hat{p}_i \leq u\} = P\{i \in J\} \cdot P\{U_1 \leq u\} = \frac{k}{s} \cdot \frac{u}{k/s} = u$$

and if $u \geq k/s$,

$$P\{\hat{p}_i \leq u\} = P\{i \in J\} \cdot 1 + P\{i \notin J\} \cdot P\{U_2 \leq u\} = \frac{k}{s} + \left(1 - \frac{k}{s}\right) \cdot \frac{u - k/s}{1 - k/s} = u.$$

Now exactly $k$ of the $\hat{p}_i$ are less than or equal to $k/s$ by construction. The prob- ability that these are all less than or equal to $\alpha k/s$ is

$$P\left\{U_1 \leq \frac{\alpha k}{s}\right\} = \frac{\alpha k/s}{k/s} = \alpha.$$

□

As is the case for the Bonferroni method, the above single-stage procedure can be strengthened by a Holm type of improvement. Consider the stepdown procedure described in (11), where now we specifically consider

$$\alpha_i = \begin{cases} \dfrac{k\alpha}{s}, & i \leq k, \\ \dfrac{k\alpha}{s+k-i}, & i > k. \end{cases} \tag{12}$$

Of course, the $\alpha_i$ depend on $s$ and $k$, but we suppress this dependence in the notation.

THEOREM 2.2. *For testing $H_i : P \in \omega_i$, $i = 1, \ldots, s$, suppose $\hat{p}_i$ satisfies* (9). *The stepdown procedure described in* (11) *with $\alpha_i$ given by* (12) *controls the $k$-FWER, that is,* (3) *holds.*

PROOF. Fix any $P$ and let $I(P)$ be the indices of the true null hypotheses. Assume $|I(P)| \geq k$ or there is nothing to prove. Order the $p$-values corresponding to the $|I(P)|$ true null hypotheses; call them

$$\hat{q}_{(1)} \leq \cdots \leq \hat{q}_{|I(P)|}.$$

Let $j$ be the smallest (random) index satisfying $\hat{p}_{(j)} = \hat{q}_{(k)}$, so

$$k \leq j \leq s - |I(P)| + k \tag{13}$$

because the largest possible index $j$ occurs when all the smallest $p$-values correspond to the $s - |I(P)|$ false null hypotheses and the next $|I(P)|$ $p$-values correspond to the true null hypotheses. So $\hat{p}_{(j)} = \hat{q}_{(k)}$. Then our generalized Holm procedure commits at least $k$ false rejections if and only if

$$\hat{p}_{(1)} \leq \alpha_1, \qquad \hat{p}_{(2)} \leq \alpha_2, \qquad \ldots, \qquad \hat{p}_{(j)} \leq \alpha_j,$$

which certainly implies that

$$\hat{q}_{(k)} = \hat{p}_{(j)} \leq \alpha_j = \frac{k\alpha}{s+k-j}.$$

But by (13),

$$\frac{k\alpha}{s+k-j} \leq \frac{k\alpha}{|I(P)|}.$$

So the probability of at least $k$ false rejections is bounded above by

$$P\left\{\hat{q}_{(k)} \leq \frac{k\alpha}{|I(P)|}\right\}.$$

By Theorem 2.1(i) the chance that the $k$th largest among $I(P)$ $p$-values is less than or equal to $k\alpha/|I(P)|$ is less than or equal to $\alpha$.　□

REMARK 2.1. Evidently, one can always reject the hypotheses corresponding to the smallest $k-1$ $p$-values without violating control of the $k$-FWER. However, it seems counterintuitive to consider a stepdown procedure whose corresponding $\alpha_i$ are not monotone nondecreasing. In addition, automatic rejection of $k-1$ hypotheses, regardless of the data, appears at the very least a little too optimistic. To ensure monotonicity, our stepdown procedure uses $\alpha_i = k\alpha/s$. Even if we were to adopt the more optimistic strategy of always rejecting the hypotheses corresponding to the first $k-1$ hypotheses, we could still only reject $k$ or more hypotheses if $\hat{p}_{(k)} \leq k\alpha/s$, which is also true for the specific procedure of Theorem 2.2.

REMARK 2.2. If the $p$-values have discrete distributions, it is possible that there may be ties among them. However, the proof remains valid regardless of how tied $p$-values are ordered because monotonicity of the $\alpha_i$ ensures that all hypotheses with a common tied $p$-value will be rejected if any of them are rejected.

The question naturally arises whether it is possible to improve the procedure further by increasing the critical values $\alpha_1, \alpha_2, \ldots$ without violating control of the $k$-FWER (3). By the previous remark we can always increase $\alpha_i$ to 1 for $i < k$. A more interesting question is whether we can increase $\alpha_i$ for $i \geq k$. We will show that this is not possible by exhibiting for each $i \geq k$ a joint distribution of the $p$-values for which

(14) $\quad P\{\hat{p}_{(1)} \leq \alpha_1, \hat{p}_{(2)} \leq \alpha_2, \ldots, \hat{p}_{(i-1)} \leq \alpha_{i-1}, \hat{p}_{(i)} \leq \alpha_i\} = \alpha.$

Moreover, changing $\alpha_i$ to $\beta_i > \alpha_i$ results in the right-hand side being greater than $\alpha$. Thus, with $i \geq k$, one cannot increase $\alpha_i$ without violating the $k$-FWER. Then, having picked $\alpha_1, \ldots, \alpha_k, \ldots, \alpha_{i-1}$, the largest possible choice for $\alpha_i$ is as stated in the algorithm.

THEOREM 2.3. (i) *Let the $\alpha_i$ be given in* (12). *For any $i \geq k$ there exists a joint distribution for $\hat{p}_1, \ldots, \hat{p}_s$ such that $s + k - i$ of the $\hat{p}_i$ are uniformly distributed on $(0,1)$ and* (14) *holds.*

(ii) *For testing $H_i : P \in \omega_i$, $i = 1, \ldots, s$, suppose $\hat{p}_i$ satisfies* (9). *For the stepdown procedure* (11) *with $\alpha_i$ given in* (12), *one cannot increase even one of the constants $\alpha_i$ (for $i \geq k$) without violating the $k$-FWER.*

Before proving the theorem, we make use of the following lemma.

LEMMA 2.1. *Fix $k$, $u$ and constants $0 < \beta_1 \leq \beta_2 \leq \cdots \leq \beta_k \leq u$. Assume for every $j = 2, \ldots, k$,*

(15) $$\frac{j(\beta_j - \beta_{j-1})}{\beta_j} \leq 1.$$

*Then there exists a joint distribution for* $(\hat{q}_1, \ldots, \hat{q}_k)$ *satisfying the* $\hat{q}_i$ *are marginally uniform on* $(0, u)$ *such that the ordered values* $\hat{q}_{(1)} \leq \cdots \leq \hat{q}_{(k)}$ *satisfy*

$$P\{\hat{q}_{(1)} \leq \beta_1, \ldots, \hat{q}_{(k)} \leq \beta_k\} = \beta_k/u. \quad (16)$$

PROOF. The proof is by induction on $k$. The result clearly holds for $k = 1$. With probability $\beta_k/u$ we will construct $(\hat{q}_1, \ldots, \hat{q}_k)$ equal to $(\tilde{q}_1, \ldots, \tilde{q}_k)$, where $\tilde{q}_i \sim U(0, \beta_k)$ for $i = 1, \ldots, k$ and such that their ordered values $\tilde{q}_{(1)} \leq \cdots \leq \tilde{q}_{(k)}$ satisfy

$$P\{\tilde{q}_{(1)} \leq \beta_1, \ldots, \tilde{q}_{(k)} \leq \beta_k\} = 1. \quad (17)$$

But, with probability $1 - \beta_k/u$, construct the $\tilde{q}_j$ to be conditionally distributed as $U(\beta_k, u)$. Then unconditionally the $\hat{q}_j$ satisfy (16) and are marginally distributed as $U(0, u)$. So it suffices to construct the $\tilde{q}_j$ satisfying $\tilde{q}_j \sim U(0, \beta_k)$ and (17).

Let $\beta_0 = 0$ and for $i = 1, \ldots, k$ let $E_i = \{(\beta_{i-1}, \beta_i]\}$ and $p_i = \beta_i - \beta_{i-1}$. First construct $Y_1, \ldots, Y_{k-1}$, each taking values in $(0, \beta_{k-1}]$ such that their ordered values $Y_{(1)} \leq \cdots \leq Y_{(k-1)}$ satisfy

$$P\{Y_{(1)} \leq \beta_1, \ldots, Y_{(k-1)} \leq \beta_{k-1}\} = 1 \quad (18)$$

and $Y_i$ is uniform on $(0, \beta_{k-1}]$. This is possible by the inductive hypothesis, since we can assume the result holds for $k-1$ as long as $\beta_1, \ldots, \beta_k$ and $u$ satisfy the stated conditions; in particular, we apply the result with $u = \beta_{k-1}$. Next, let $Y_k$ be uniform on $E_i$ with probability $\theta p_i$ for $i = 1, \ldots, k-1$ and let it be uniform on $E_k$ with probability $1 - \theta\beta_{k-1}$, where $\theta$ satisfies

$$\theta = \frac{1}{\beta_{k-1}}\left[1 - \frac{k(\beta_k - \beta_{k-1})}{\beta_k}\right]. \quad (19)$$

Finally, let $\tilde{q}_1, \ldots, \tilde{q}_k$ be a random permutation of $Y_1, \ldots, Y_k$. Because of (18) and the fact that $Y_k \leq \beta_k$, the ordered values of $Y_1, \ldots, Y_k$ and hence the ordered values of $\tilde{q}_1, \ldots, \tilde{q}_k$ satisfy (17). Furthermore, it is easy to check that $\tilde{q}_i$ falls in $E_j$ with probability $p_j$ and so $\tilde{q}_i$ is $U(0, \beta_k)$. Indeed, if $j < k$, the probability that $\tilde{q}_i$ falls in $E_j$, conditional on $\tilde{q}_i$ not being equal to $Y_k$, is $p_i/\beta_{k-1}$ and is $\theta p_i$ in the latter case, which unconditionally is

$$\frac{k-1}{k} \cdot \frac{p_i}{\beta_{k-1}} + \frac{1}{k}\theta p_i = p_i,$$

and similarly for the probability that $\hat{q}_i$ falls in $E_k$. The only detail that remains is to note that this construction with $\theta$ defined in (19) is possible only if $\theta p_i$ and $1 - \theta\beta_{k-1}$ are all values in $(0, 1)$. But

$$1 - \theta\beta_{k-1} = \frac{k(\beta_k - \beta_{k-1})}{\beta_k},$$

which is certainly $\geq 0$ since $\beta_k \geq \beta_{k-1}$. It is also $\leq 1$ by the assumption (15). Also,

$$\theta p_i = \frac{p_i}{\beta_{k-1}} \cdot \left[1 - \frac{k(\beta_k - \beta_{k-1})}{\beta_k}\right].$$

But the first factor $p_i/\beta_{k-1}$ is in $(0,1)$ as is the latter by the above, and so the product is in $(0,1)$. □

PROOF OF THEOREM 2.3. The case $i = k$ follows from the construction in the proof of Theorem 2.1. Let the first $i-k$ of the $\hat{p}_j$ be identically equal to 0. (Actually, rather than point mass at 0, any distribution supported on $[0, \alpha_1)$ will do.) For the remaining $s' = s+k-i$ $p$-values $\hat{p}_j$, $j = i-k+1, \ldots, s$, randomly choose $k$ indices from $i-k+1, \ldots, s$. The $k$ that are chosen will be marginally $U(0, k/s')$ and have a joint distribution which will be specified below; the remaining $s-i$ can be taken to be distributed as $U(k/s', 1)$.

Let $\hat{q}_1, \ldots, \hat{q}_k$ denote the $k$ observations that are marginally $U(0, k/s')$. We need to specify the joint distribution of $\hat{q}_1, \ldots, \hat{q}_k$ so that their ordered values $\hat{q}_{(1)} \leq \cdots \leq \hat{q}_{(k)}$ satisfy

$$(20) \qquad P\{\hat{q}_{(1)} \leq \alpha_{i-k+1}, \hat{q}_{(2)} \leq \alpha_{i-k+2}, \ldots, \hat{q}_{(k)} \leq \alpha_i\} = \alpha$$

(because $\hat{q}_{(j)} = \hat{p}_{(j+i-k)}$ for $j = 1, \ldots, k$). So the problem reduces to constructing a joint distribution for $(\hat{q}_1, \ldots, \hat{q}_k)$ satisfying (20) subject to the constraint that $\hat{q}_j$ is marginally distributed as $U(0, k/s')$. To do this, apply Lemma 2.1 with $u = k/s'$ and $\beta_j = \alpha_{i-k+j}$. We need to verify the conditions of the lemma, which reduces to showing

$$(21) \qquad \frac{j(\alpha_{i-k+j} - \alpha_{i-k+j-1})}{\alpha_{i-k+j}} \leq 1$$

for $i \geq k$ (and $s$ and $k$ fixed). But, if $i - k + j - 1 \leq k$, then the left-hand side of (21) is 0; otherwise it is easily seen to simplify to

$$(22) \qquad \frac{j}{s+2k-i-j} \leq \frac{j}{s+k-j} \leq k/s,$$

where the first inequality holds because $i \geq k$ and the second because $j \leq k$. But $k/s \leq 1$ and so the conditions of the lemma are satisfied. Therefore, we can conclude that the left-hand side of (20) is given by

$$\frac{\beta_k}{u} = \frac{\alpha_i}{k/s'} = \alpha,$$

and (i) is proved.

To prove (ii), the construction used in (i) can be used even if $\alpha_i$ is replaced by $\bar{\alpha}_i > \alpha_i$, as long as such a switch still allows one to appeal to the lemma. However, the same argument works as long as $\bar{\alpha}_i$ does not get bigger than $s/k \cdot \alpha_i$, so that the argument leading to (22) being less than or equal to 1 still applies. For such an $\bar{\alpha}_i$, the argument for (i) then shows that, if the left-hand side of (14) has $\alpha_i$ replaced by $c\alpha_i$ for some $1 < c < s/k$, then the right-hand side of (14) will be $c\alpha > \alpha$, which would violate control of the $k$-FWER. □

**3. Control of the false discovery proportion.** The number $k$ of false rejections that one is willing to tolerate will often increase with the number of hypotheses rejected. So it might be of interest to control not the number of false rejections (sometimes called false discoveries) but the proportion of false discoveries. Specifically, let the *false discovery proportion* (FDP) be defined by

$$(23) \quad FDP = \begin{cases} \dfrac{\text{Number of false rejections}}{\text{Total number of rejections}}, & \text{if the denominator} \\ & \text{is greater than 0}, \\ 0, & \text{if there are no rejections.} \end{cases}$$

Thus FDP is the proportion of rejected hypotheses that are rejected erroneously. When none of the hypotheses is rejected, both numerator and denominator of that proportion are 0; since in particular there are no false rejections, the FDP is then defined to be 0.

Benjamini and Hochberg [1] proposed to replace control of the FWER by control of the *false discovery rate* (FDR), defined as

$$(24) \quad FDR = E(FDP).$$

The FDR has gained wide acceptance in both theory and practice, largely because Benjamini and Hochberg proposed a simple stepup procedure to control the FDR. Unlike control of the $k$-FWER, however, their procedure is not valid without assumptions on the dependence structure of the $p$-values. Their original paper assumed the very strong assumption of independence of $p$-values, but this has been weakened to include certain types of dependence; see [2]. In any case, control of the FDR does not prohibit the FDP from varying, even if its average value is bounded. Instead, we consider an alternative measure of control that guarantees the FDP is bounded, at least with prescribed probability. That is, for a given $\gamma$ and $\alpha$ in $(0,1)$, we require

$$(25) \quad P\{FDP > \gamma\} \leq \alpha.$$

To develop a stepdown procedure satisfying (25), let $F$ denote the number of false rejections. At step $i$, having rejected $i-1$ hypotheses, we want to guarantee $F/i \leq \gamma$, that is, $F \leq \lfloor \gamma i \rfloor$, where $\lfloor x \rfloor$ is the greatest integer less than or equal to $x$. So, if $k = \lfloor \gamma i \rfloor + 1$, then $F \geq k$ should have probability no greater than $\alpha$; that is, we must control the number of false rejections to be less than or equal to $k$. Therefore, we use the stepdown constant $\alpha_i$ with this choice of $k$ (which now depends on $i$); that is,

$$\alpha_i = \frac{(\lfloor \gamma i \rfloor + 1)\alpha}{s + \lfloor \gamma i \rfloor + 1 - i}. \tag{26}$$

We give two results that show the stepdown procedure with this choice of $\alpha_i$ satisfies (25). Unfortunately, like FDR control, some assumptions on the dependence of $p$-values are required, at least by our method of proof. Later, we will modify the method so we can dispense with the dependence assumptions. As before, $\hat{p}_1, \ldots, \hat{p}_s$ denotes the $p$-values of the individual tests. Also, let $\hat{q}_1, \ldots, \hat{q}_{|I|}$ denote the $p$-values corresponding to the $|I| = |I(P)|$ true null hypotheses. So $q_i = p_{j_i}$, where $j_1, \ldots, j_{|I|}$ correspond to the indices of the true null hypotheses. Also, let $\hat{r}_1, \ldots, \hat{r}_{s-|I|}$ denote the $p$-values of the false null hypotheses. Consider the following condition: for any $i = 1, \ldots, |I|$,

$$P\{\hat{q}_i \leq u | \hat{r}_1, \ldots, \hat{r}_{s-|I|}\} \leq u; \tag{27}$$

that is, conditional on the observed $p$-values of the false null hypotheses, a $p$-value corresponding to a true null hypothesis is (conditionally) dominated by the uniform distribution, as it is unconditionally in the sense of (7). No assumption is made regarding the unconditional (or conditional) dependence structure of the true $p$-values, nor is there made any explicit assumption regarding the joint structure of the $p$-values corresponding to false hypotheses, other than the basic assumption (27). So, for example, if the $p$-values corresponding to true null hypotheses are independent of the false ones, but have arbitrary joint dependence within the group of true null hypotheses, the above assumption holds.

THEOREM 3.1. *Assume condition* (27). *Then the stepdown procedure with $\alpha_i$ given by* (26) *controls the FDP in the sense of* (25).

PROOF. Assume the number of true null hypotheses is $|I(P)| > 0$ (or there is nothing to prove) and the number of false null hypotheses is $f = s - |I(P)|$. The argument is conditional on the $\{\hat{r}_i\}$. Let

$$\hat{r}_{(1)} \leq \hat{r}_{(2)} \leq \cdots \leq \hat{r}_{(f)}$$

denote the ordered values of the $\hat{r}_i$ and similarly for the $\hat{q}_i$. Let $\alpha_0 = 0$ and define $R_i$ to be the number of $\hat{r}_i$ in the interval $(\alpha_{i-1}, \alpha_i]$. (Actually, assume $R_1$ includes the value 0 as well.) Given the values of $\hat{r}_1, \ldots, \hat{r}_f$, it may be impossible to have $FDP > \gamma$, that is,

$$P\{FDP > \gamma | \hat{r}_1, \ldots, \hat{r}_f\} = 0.$$

Otherwise, let $j = j(\hat{r}_1, \ldots, \hat{r}_f)$ be defined as

$$j = \min\left\{m : m - \sum_{i=1}^{m} R_i > m\gamma\right\}. \tag{28}$$

To interpret this, given the $p$-values of the false hypotheses, $j$ is the smallest critical index (depending only on the $\hat{r}_i$) where it is possible to have $FDP > \gamma$, except whenever there are several $p$-values within an interval $(\alpha_{i-1}, \alpha_i)$ we consider the index of the largest one. The point of the construction is that if the stepdown procedure stops at an index $m < j$, then $m - \sum_i R_i/m \leq \gamma$ and so $FDP \leq \gamma$. On the other hand, if the event $FDP > \gamma$ occurs, then there must be a rejection of a true null hypothesis at step $j$.

For example, if $s = 100$, $f = 5$ and $\gamma = 0.1$, then if all five of the $\hat{r}_i$ are less than $\alpha_1$, then we define $j = 6$ even though the smallest true $p$-value could be the smallest among the 100. So the FDP could be greater than 0.1 after the first step of the algorithm if $\hat{q}_{(1)} < \hat{r}_{(1)}$, but even if this is the case, we then know we will reject at least six total hypotheses. So the important point here is that, given such a configuration of $\{\hat{r}_i\}$, in order for FDP to be greater than 0.1, it must be the case that we reject a true null hypothesis at step 6.

Note that, with $j$ so defined, $R_j = 0$. For if $\sum_{i=1}^{j} R_i = j - k$ with $k/j > \gamma$ and $R_j > 0$, then

$$\sum_{i=1}^{j-1} R_i = j - k - R_j \leq j - 1 - k$$

and $k/(j-1) > \gamma$, so that $m = j - 1$ satisfies the criterion. Furthermore, we also have $\sum_{i=1}^{j} R_i = j - k$ (so not $< j - k$), where $k/j > \gamma$, because if $\sum_{i=1}^{j} R_i < j - k \leq j - 1 - k$ say, then $k/(j-1) > \gamma$ if $k/j > \gamma$ and so $j$ can again be reduced to $j - 1$.

In addition, at the index $j$ it must be the case that

$$k = k(j) = j - \sum_{i=1}^{j} R_i = 1 + \lfloor \gamma j \rfloor.$$

But $k > \gamma j$ implies $k \geq \lfloor \gamma j \rfloor + 1$. But if $k > \lfloor \gamma j \rfloor + 1$, then $k - 1 \geq \lfloor \gamma j \rfloor + 1$ and so

$$\frac{k-1}{j-1} \geq \frac{\lfloor \gamma j \rfloor + 1}{j-1} > \gamma,$$

the last equality trivially following from $1 + \lfloor \gamma j \rfloor \geq \gamma j > \gamma(j-1)$.

We can now complete the argument. At the index $j$ we must have $k = j - \sum_{i=1}^{j} R_i = 1 + \lfloor \gamma j \rfloor$ of the $\hat{q}_i$ being $\leq \alpha_j$. But from Theorem 2.1 (applied conditional on the $\hat{r}_i$),

$$P\{\text{at least } k(j) \text{ of the } \hat{q}_i \leq \alpha_j | \hat{r}_1, \ldots, \hat{r}_f\}$$
$$\leq \frac{|I|\alpha_j}{k(j)}$$
$$= \frac{|I|(\lfloor \gamma j \rfloor + 1)\alpha}{k(j)(s + \lfloor \gamma j \rfloor + 1 - j)} = \frac{|I|\alpha}{s + \lfloor \gamma j \rfloor + 1 - j}.$$

But $|I| \leq s - \sum_{i=1}^{j} R_i = s - j + k$, so the above probability is less than or equal to

$$\frac{s - j + k}{s + \lfloor \gamma j \rfloor + 1 - j} \cdot \alpha = \alpha.$$

Therefore,

$$P\{FDP > \gamma | \hat{r}_1, \ldots, \hat{r}_f\} \leq \alpha,$$

which of course implies $P\{FDP > \gamma\} \leq \alpha$. □

Next, we prove the same stepdown procedure controls the FDP in the sense of (25) under an alternative assumption. Here, the assumption only involves the dependence of the $p$-values corresponding to true null hypotheses.

THEOREM 3.2. *Consider testing $s$ null hypotheses, with $|I|$ of them true. Let $\hat{q}_{(1)} \leq \cdots \leq \hat{q}_{(|I|)}$ denote their corresponding ordered p-values. Set $M = \min(\lfloor \gamma s \rfloor + 1, |I|)$.*

(i) *For the stepdown procedure with $\alpha_i$ given by (26),*

(29) $$P\{FDP > \gamma\} \leq P\left\{\bigcup_{i=1}^{M}\left\{\hat{q}_{(i)} \leq \frac{i\alpha}{|I|}\right\}\right\}.$$

(ii) *Therefore, if the joint distribution of the p-values of the true null hypotheses satisfies Simes inequality, that is,*

$$P\left\{\left\{\hat{q}_{(1)} \leq \frac{\alpha}{|I|}\right\} \cup \left\{\hat{q}_{(2)} \leq \frac{2\alpha}{|I|}\right\} \cup \cdots \cup \{\hat{q}_{(|I|)} \leq \alpha\}\right\} \leq \alpha,$$

*then $P\{FDP > \gamma\} \leq \alpha$.*

PROOF. Let $j$ be the smallest (random) index where the FDP exceeds $\gamma$ for the first time at step $j$; that is, the number of false rejections corresponding to the first $j-1$ rejections divided by $j$ exceeds $\gamma$ for the first time at $j$. If $j$ is such that $\gamma j < 1$, then $FDP > \gamma$ at step $j$ implies $\hat{p}_{(j)} \leq \alpha_j$. But this implies

$$\hat{q}_{(1)} \leq \alpha_j = \frac{\alpha}{s+1-j} \leq \frac{\alpha}{|I|},$$

because the number of true null hypotheses $|I|$ necessarily satisfies $|I| \leq s - (j-1)$ for such a $j$.

Similarly, if $j$ is such that $1 \leq \gamma j < 2$, then we must have $\hat{p}_{(i)} \leq \alpha_i$ and $\hat{p}_{(j)} \leq \alpha_j$ for some $i < j$, where $i, j$ correspond to true null hypotheses. But for such a $j$, $\alpha_j = 2\alpha/(s+2-j)$, and so we must have $\hat{q}_{(2)} \leq 2\alpha/(s-j+2)$. But, by definition of $j$, we must have $|I| \leq s - (j-2)$ and so $\hat{q}_{(2)} \leq 2\alpha/|I|$.

Continuing in this way, if $m - 1 \leq \gamma j < m$, the event $FDP > \gamma$ at step $j$ implies $\hat{q}_{(m)} \leq m\alpha/|I|$. The largest value of $j$ is of course $s$ and so the largest possible $m$ is $\lfloor \gamma s \rfloor + 1$. Also, we cannot have $m > |I|$. So, with $M$ as in the statement of the theorem,

$$P\{FDP > \gamma\} \leq \sum_{m=1}^{M} P\left\{\hat{q}_{(m)} \leq \frac{m\alpha}{|I|}, m-1 \leq \gamma j < m\right\}$$

$$\leq \sum_{m=1}^{M} P\left\{\bigcup_{i=1}^{M}\left\{\hat{q}_{(i)} \leq \frac{i\alpha}{|I|}\right\}, m-1 \leq \gamma j < m\right\}$$

$$\leq P\left\{\bigcup_{i=1}^{M}\left\{\hat{q}_{(i)} \leq \frac{i\alpha}{|I|}\right\}\right\}.$$

Part (ii) follows trivially. $\square$

In fact, there are many joint distributions of positively dependent variables for which Simes inequality is known to hold. In particular, Sarkar and Chang [11] and Sarkar [10] have shown that the Simes inequality holds for the family of distributions which is characterized by the multivariate positive of order 2 condition, as well as some other important distributions.

Theorem 3.2 points toward a method that controls the FDP without any dependence assumptions. One simply needs to bound the right-hand side of (29). In fact, Hommel [6] has shown that

$$P\left\{\bigcup_{i=1}^{|I|}\left\{\hat{q}_{(i)} \leq \frac{i\alpha}{|I|}\right\}\right\} \leq \alpha \sum_{i=1}^{|I|} \frac{1}{i}.$$

This suggests we replace $\alpha$ by $\alpha(\sum_{i=1}^{|I|}(1/i))^{-1}$. But of course $|I|$ is unknown. So one possibility is to bound $|I|$ by $s$, which then results in replacing $\alpha$ by $\alpha/C_s$, where

(30) $$C_j = \sum_{i=1}^{j}(1/i).$$

As is well known, $C_s \approx \log(s + 0.5) + \zeta_E$, with $\zeta_E \approx 0.5772156649$ known as Euler's constant. Clearly, changing $\alpha$ in this way is much too conservative and results in a much less powerful method. However, notice in (29) that we really only need to bound the union over $M \leq \lfloor \gamma s + 1 \rfloor$ events. Therefore, we need to slightly generalize the inequality by Hommel [6], which is done in the following lemma.

LEMMA 3.1. *Suppose $\hat{p}_1, \ldots, \hat{p}_t$ are p-values in the sense that $P\{\hat{p}_i \leq u\} \leq u$ for all $i$ and $u$ in $(0,1)$. Let their ordered values be $\hat{p}_{(1)} \leq \cdots \leq \hat{p}_{(t)}$. Let $0 = \beta_0 \leq \beta_1 \leq \beta_2 \leq \cdots \leq \beta_m \leq 1$ for some $m \leq t$.*

(i) *Then*

$$\text{(31)} \quad P\{\{\hat{p}_{(1)} \leq \beta_1\} \cup \{\hat{p}_{(2)} \leq \beta_2\} \cup \cdots \cup \{\hat{p}_{(m)} \leq \beta_m\}\} \leq t \sum_{i=1}^{m} (\beta_i - \beta_{i-1})/i.$$

(ii) *As long as the right-hand side of (31) is less than or equal to 1, the bound is sharp in the sense that there exists a joint distribution for the p-values for which the inequality is an equality.*

PROOF. Let $J$ be the smallest (random) index $j$ among $1 \leq j \leq m$ for which $\hat{p}_{(j)} \leq \beta_j$; define $J$ to be $t+1$ if $\hat{p}_{(j)} > \beta_j$ for all $1 \leq j \leq m$. Let $\theta_k = P\{J = k\}$. Then the left-hand side of (31) is equal to

$$P\left\{\bigcup_{k=1}^{m} \{J = k\}\right\} = \sum_{k=1}^{m} \theta_k,$$

since the events $\{J = k\}$ are disjoint. We wish to bound $\sum_k \theta_k$. For any $1 \leq j \leq m$,

$$\sum_{k=1}^{j} JI\{J = k\} = JI\{J \leq j\} \leq S_j,$$

where $S_j$ is the number of p-values $\leq \beta_j$. Taking expectations yields

$$\text{(32)} \quad \sum_{k=1}^{j} k\theta_k \leq t\beta_j, \qquad j = 1, \ldots, m.$$

For $j = 1, \ldots, m-1$, multiply both sides of (32) by $1/[j(j+1)]$, and for $j = m$, multiply both sides by $1/m$; then sum over $j$ to yield

$$\text{(33)} \quad \sum_{j=1}^{m-1} \frac{1}{j(j+1)} \sum_{k=1}^{j} k\theta_k + \frac{1}{m} \sum_{k=1}^{m} k\theta_k \leq \sum_{j=1}^{m-1} \frac{t\beta_j}{j(j+1)} + \frac{t\beta_m}{m}.$$

By changing the order of summation, the left-hand side of (33) becomes

$$\sum_{k=1}^{m-1} k\theta_k \left(\frac{1}{k} - \frac{1}{m}\right) + \frac{1}{m}\sum_{k=1}^{m} k\theta_k = \sum_{k=1}^{m} \theta_k.$$

The right-hand side of (33) is easily seen to be the right-hand side of (31) and (i) follows.

To prove (ii), we construct $\hat{p}_1, \ldots, \hat{p}_t$ as follows. Let $U_i$ be uniform in $I_i$ and let $U_{m+1}$ be uniform in $(\beta_m, 1)$. Let $p$ be equal to the right-hand side of (31), assumed less than or equal to 1. Let $\pi_1, \ldots, \pi_m$ be probabilities summing to 1, with $\pi_i \propto (\beta_i - \beta_{i-1})/i$. Then, with probability $\pi_i p$, randomly pick $i$ indices and let those $p$-values be equal to $U_i$, and the remaining $t - i$ $p$-values equal to $U_{m+1}$. With the remaining probability $1 - p$, let all $p$-values be equal to $U_{m+1}$. With this construction it is easily checked that $\hat{p}_i$ is uniform on $(0, 1)$ and the left-hand side of (31) is equal to the right-hand side of (31). □

Theorem 3.2 and Lemma 3.1 now lead to the following result.

THEOREM 3.3. *For testing $H_i : P \in \omega_i$, $i = 1, \ldots, s$, suppose $\hat{p}_i$ satisfies (9). Consider the stepdown procedure with constants $\alpha'_i = \alpha_i / C_{(\lfloor \gamma s \rfloor + 1)}$, where $\alpha_i$ is given by (26) and $C_j$ is defined by (30). Then $P\{FDP > \gamma\} \leq \alpha$.*

PROOF. By Theorem 3.2(i), $P\{FDP > \gamma\}$ is bounded by the right-hand side of (29) with $\alpha$ replaced by $\alpha/C_{\lfloor \gamma s \rfloor + 1}$, which is further bounded by the same expression with $M$ replaced by $\lfloor \gamma s \rfloor + 1$. Then apply Lemma 3.1 with $t = |I|$ and $\beta_i = i\alpha/(C_{\lfloor \gamma s \rfloor + 1}|I|)$. □

It is of interest to compare control of the FDP with control of the FDR. Some obvious connections between methods that control the FDP in the sense of (25) and methods that control its expected value, the FDR, can be made. Indeed, for any random variable $X$ on $[0, 1]$, we have

$$E(X) = E(X|X \leq \gamma)P\{X \leq \gamma\} + E(X|X > \gamma)P\{X > \gamma\}$$
$$\leq \gamma P\{X \leq \gamma\} + P\{X > \gamma\},$$

which leads to

(34) $$\frac{E(X) - \gamma}{1 - \gamma} \leq P\{X > \gamma\} \leq \frac{E(X)}{\gamma},$$

with the last inequality just Markov's inequality. Applying this to $X = FDP$, we see that, if a method controls the FDR at level $q$, then it controls the FDP in the sense $P\{FDP > \gamma\} \leq q/\gamma$. Obviously, this is very crude because if $q$ and $\gamma$ are both small, the ratio can be quite large. The first inequality in (34) says that if the FDP is controlled in the sense of (25), then the FDR is controlled at level $\alpha(1-\gamma)+\gamma$, which is greater than or equal to $\alpha$ but typically only slightly. These crude arguments suggest that control of the FDP is perhaps more stringent than control of the FDR.

The comparison of actual methods, however, is complicated by the fact that the FDR controlling procedure of Benjamini and Hochberg [1] is a stepup procedure, but we have only considered stepdown procedures. It is interesting to note that, in order to make our procedure work without any dependence assumptions, we needed to change $\alpha$ to $\alpha/C_{\lfloor \gamma s \rfloor + 1}$. Benjamini and Yekutieli [2] show that the Benjamini–Hochberg procedure that controls the FDR at level $q$ can also work without dependence assumptions, if you replace $q$ by $q/C_s$. Clearly, this is a more drastic change since $C_s$ is typically much larger than $C_{\lfloor \gamma s \rfloor + 1}$. Such connections need to be explored more fully.

**4. Conclusions.** We have seen that a very simple stepdown procedure is available to control the $k$-FWER under absolutely no assumptions on the dependence structure of the $p$-values. Furthermore, control of the $k$-FWER provides a measure of control for the *actual* number of false rejections, while the number of false rejections in the case of the FDR can vary widely. We have also considered two stepdown methods that control the FDP in the sense of (25). The first method provides control under very reasonable types of dependence assumptions, while the second holds in general.

**Acknowledgments.** We thank Juliet Shaffer and Michael Wolf for some helpful discussion and references. We also thank the referees and an Associate Editor for many helpful suggestions that greatly improved the clarity of the paper. Thanks to Wenge Guo for pointing out an error in an earlier version.

After the revision and acceptance of this paper, we became aware of the work by Hommel and Hoffman [7] which has much overlap with the results in Section 2, and we'd like to thank Helmut Finner for pointing out this oversight. In particular, Hommel and Hoffman [7] provide Theorem 2.1(i) with proof, Theorem 2.2 (stated but no proof) and a weaker version of Theorem 2.3(ii) (stated but no proof). They attribute the idea of controlling the number of false hypotheses to Victor [14], who also suggested control of the FDP. However, Hommel and Hoffman did not further discuss control of the FDP as they "could not find suitable procedures satisfying this criterion." As far as we know, the three theorems in Section 3 which address control of the FDP are new.

DEPARTMENT OF STATISTICS
UNIVERSITY OF CALIFORNIA
BERKELEY, CALIFORNIA 94720
USA

DEPARTMENT OF STATISTICS
STANFORD UNIVERSITY
STANFORD, CALIFORNIA 94305-4065
USA
E-MAIL: romano@stat.stanford.edu